\documentclass[11pt]{amsart}
\usepackage{amscd}
\usepackage{amsmath}
\usepackage{amsxtra}
\usepackage{amsfonts}
\usepackage{amssymb}

\newtheorem{theorem}{Theorem}[section]
\newtheorem{corollary}[theorem]{Corollary}
\newtheorem{lemma}[theorem]{Lemma}
\newtheorem{proposition}[theorem]{Proposition}

\theoremstyle{definition}
\newtheorem{definition}[theorem]{Definition}
\newtheorem{remark}[theorem]{Remark}

\theoremstyle{remark}

\renewcommand{\theclaim}{\textup{\theclaim}}

\newtheorem*{acknowledgements}{Acknowledgements}

\numberwithin{equation}{section}

\def\openone

{\mathchoice

{\hbox{\upshape \small1\kern-3.3pt\normalsize1}}

{\hbox{\upshape \small1\kern-3.3pt\normalsize1}}

{\hbox{\upshape \tiny1\kern-2.3pt\SMALL1}}

{\hbox{\upshape \Tiny1\kern-2pt\tiny1}}}

\makeatletter

\newbox\ipbox

\newcommand{\ip}[2]{\left\langle #1\, , \,#2\right\rangle}
\newcommand{\diracb}[1]{\left\langle #1\mathrel{\mathchoice

{\setbox\ipbox=\hbox{$\displaystyle \left\langle\mathstrut
#1\right.$}

\vrule height\ht\ipbox width0.25pt depth\dp\ipbox}

{\setbox\ipbox=\hbox{$\textstyle \left\langle\mathstrut
#1\right.$}

\vrule height\ht\ipbox width0.25pt depth\dp\ipbox}

{\setbox\ipbox=\hbox{$\scriptstyle \left\langle\mathstrut
#1\right.$}

\vrule height\ht\ipbox width0.25pt depth\dp\ipbox}

{\setbox\ipbox=\hbox{$\scriptscriptstyle \left\langle\mathstrut
#1\right.$}

\vrule height\ht\ipbox width0.25pt depth\dp\ipbox}

}\right. }

\newcommand{\dirack}[1]{\left. \mathrel{\mathchoice

{\setbox\ipbox=\hbox{$\displaystyle \left.\mathstrut
#1\right\rangle$}

\vrule height\ht\ipbox width0.25pt depth\dp\ipbox}

{\setbox\ipbox=\hbox{$\textstyle \left.\mathstrut
#1\right\rangle$}

\vrule height\ht\ipbox width0.25pt depth\dp\ipbox}

{\setbox\ipbox=\hbox{$\scriptstyle \left.\mathstrut
#1\right\rangle$}

\vrule height\ht\ipbox width0.25pt depth\dp\ipbox}

{\setbox\ipbox=\hbox{$\scriptscriptstyle \left.\mathstrut
#1\right\rangle$}

\vrule height\ht\ipbox width0.25pt depth\dp\ipbox}

} #1\right\rangle}

\newcommand{\cj}[1]{\overline{#1}}

\newcommand{\bz}{\mathbb{Z}}

\newcommand{\br}{\mathbb{R}}
\newcommand{\bc}{\mathbb{C}}
\newcommand{\bt}{\mathbb{T}}
\newcommand{\bn}{\mathbb{N}}

\def\blfootnote{\xdef\@thefnmark{}\@footnotetext}

\renewcommand{\mod}{\operatorname{mod}}

\input xy
\xyoption{all}
\usepackage{amssymb}

\newcommand{\Span}{\overline{\operatorname*{span}}}

\def\be{\mathbb{E}}
\def\H{\mathcal{H}}

\def\-{^{-1}}

\begin{document}

\title[Irreducible wavelet representations]{Irreducible wavelet representations and ergodic automorphisms on solenoids}
\author{Dorin Ervin Dutkay}
\blfootnote{This research was supported in part by the Swedish Foundation for International Cooperation in Research and Higher Education (STINT) and the Swedish Research Council. The first author was also supported by a NSF Young Investigator Award, at Texas A\&M. }
\address{[Dorin Ervin Dutkay] University of Central Florida\\
	Department of Mathematics\\
	4000 Central Florida Blvd.\\
	P.O. Box 161364\\
	Orlando, FL 32816-1364\\
U.S.A.\\} \email{ddutkay@mail.ucf.edu}
\author{David R. Larson}
\address{[David R. Larson] Texas A\& M University\\
Department of Mathematics\\
College Station, TX 77832\\
U.S.A.\\}
\email{larson@math.tamu.edu}
\author{Sergei Silvestrov}
\address{[Sergei Silvestrov] Centre for Mathematical Sciences\\ 
Lund University\\
Box 118, SE-221 00 Lund, Sweden}
\email{Sergei.Silvestrov@math.lth.se}
\thanks{} 
\subjclass[2000]{42C40 ,28D05,47A67,28A80 }
\keywords{Representation, ergodic automorphism, Cantor set, Ruelle operator, solenoid, refinable function}

\begin{abstract}
We focus on the irreducibility of wavelet representations. We present some connections between the following notions: covariant wavelet representations, ergodic shifts on solenoids, fixed points of transfer (Ruelle) operators and solutions of refinement equations. We investigate the irreducibility of the wavelet representations, in particular the representation associated to the Cantor set, introduced in \cite{DuJo06b}, and we present several equivalent formulations of the problem. 
\end{abstract}
\maketitle \tableofcontents
\section{Introduction}
The interplay between dynamical and systems and operator theory is now a well developed subject \cite{Tom92, Fur99,BrJo91,Con94}. In particular, the operator theoretic approach to wavelet theory has been extremely productive \cite{DaLa98,HaLa00,BrJo99,BrJo}. We will work along the same lines: we are interested in the connections between irreducible covariant representations, ergodic shifts on solenoids and fixed points of transfer (or Ruelle) operators.

\subsection{Classical wavelet theory}
In the theory of wavelets (see e.g., \cite{Dau92}), orthonormal bases for $L^2(\br)$ are constructed by applying dilation and translation operators, in a certain order, to a given vector $\psi$ called the {\it wavelet}. Thus from the start of this construction, we have two unitary operators:
$$Uf(x)=\frac{1}{\sqrt{2}}f\left(\frac x2\right),\quad Tf(x)=f(x-1),\quad(f\in L^2(\br),x\in\br)$$
which satisfy a {\it covariance relation}:
$$UTU^{-1}=T^2.$$
Using Borel functional calculus, one can define a representation of $L^\infty(\bt)$, where $\bt$ is the unit circle:
$$\pi(f)=f(T)$$
so in particular $\pi(z^n)=T^n$, and this representation will satisfy the covariance relation
\begin{equation}\label{eqcov}
U\pi(f)U^{-1}=\pi(f(z^2)),\quad(f\in L^\infty(\bt))
\end{equation}

The main technique of constructing wavelets is by {\it multiresolutions}: one starts with a {\it quadrature-mirror-filter (QMF)} $m_0\in L^\infty(\bt)$, ($\bt$ is the unit circle)  that satisfies the {\it QMF-condition}
$$\frac{1}{2}\sum_{w^2=z}|m_0(w)|^2=1,\quad(z\in\bt),$$
the {\it low-pass condition} $m_0(1)=\sqrt{2}$, and perhaps some regularity (Lipschitz, etc.)

Then, a {\it scaling function} is constructed by an infinite product formula 
$$\hat\varphi(x)=\prod_{n=1}^\infty\frac{ m_0\left(e^{2\pi i\frac x{2^n}}\right)}{\sqrt{2}},$$
where we denote by $\hat f$ the Fourier transform of the function $f$
$$\hat f(x)=\int_{\br}f(t)e^{-2\pi itx}\,dt,\quad(x\in\br).$$
\begin{definition}\label{defi1}
We call the function $\varphi$ the {\it scaling function} associated to the QMF $m_0$. 
The scaling function satisfies the {\it scaling equation}
\begin{equation}\label{eqscaling}
U\varphi=\pi(m_0)\varphi,
\end{equation}
and it generates a sequence of subspaces $V_n$, $n\in\bz$: 
$$V_0=\Span\{T^k\varphi\,|\,k\in\bz\}=\Span\{\pi(f)\varphi\,|\,f\in L^\infty(\bt)\},$$
$$V_n=U^{-n}V_0,\quad(n\in\bz).$$
We call $(V_n)_{n\in \bz}$ the {\it multiresolution} associated to $\varphi$. The multiresolution has the properties that $V_n\subset V_{n+1}$ (this follows from the scaling equation), 
\begin{equation}
\cj{\bigcup_{n\in\bz} V_n}=L^2(\br).
\label{eqdense}
\end{equation}
If $m_0$ is carefully chosen, one gets an {\it orthonormal scaling function} $\varphi$, i.e., its translates are orthogonal 
$$\ip{T^k\varphi}{T^l\varphi}=\delta_{kl},\quad(k,l\in\bz).$$
Equivalently
\begin{equation}
\ip{\pi(f)\varphi}{\varphi}=\int_{\bt}f\,d\mu,\quad(f\in L^\infty(\bt))
\label{eqortho}
\end{equation}

Once the orthonormal scaling function and the multiresolution are constructed the wavelet is obtained by considering the {\it detail space} $W_0:=V_1\ominus V_0$. Analyzing the multiplicity of the representation $\pi$ on the spaces $V_0$ and $V_1$, one can see that there is a function $\psi$ such that $\{T^k\psi\,|\,k\in\bz\}$ is an orthonormal basis for $W_0$. Applying $U^n$, one gets that 
$$\{U^nT^k\psi\,|\,n,k\in\bz\}$$
is an orthonormal basis for $L^2(\br)$, thus $\psi$ is a {\it wavelet}.

\subsection{Wavelets on the Cantor set}\label{secwc}
Let $\mathbf C$ be the Middle Third Cantor set. A quick inspection shows that its characteristic function satisfies the following scaling equation:
$$\chi_{\mathbf C}\left(\frac{x}{3}\right)=\chi_{\mathbf C}(x)+\chi_{\mathbf C}(x-2),\quad(x\in\br).$$
This enables one to construct a multiresolution structure where $\chi_{\mathbf C}$ is a scaling function, not in $L^2(\br)$ where $\mathbf C$ has measure zero, but in $L^2$ of a Hausdorff measure (see \cite{DuJo06b}). More precisely, let 
$$\mathcal R:=\bigcup\left\{\mathbf C+\frac k{3^n}\,|\, k,n\in\bz\right\}$$
and let $\mathfrak H^s$ be the Hausdorff measure associated to the Hausdorff dimension $s=\log_32$ of the Cantor set, restricted to $\mathcal R$.

Recall (see \cite{Fal03}) that the Hausdorff measure for dimension $s$ is defined as follows: for a subset $E$ of $\br$, define for $\delta>0$:
$$\mathfrak H^s_\delta(E):=\inf\left\{\sum_{i\in I}\mbox{diam}(A_i)^s : E\subset\bigcup_{i\in I}A_i\, , \mbox{diam}(A_i)<\delta\right\}.$$
Then 
$$\mathfrak H^s(E):=\lim_{\delta\rightarrow0}\mathfrak H^s_\delta(E)$$
defines a metric outer measure. The Hausdorff measure is the restriction of $\mathfrak H^s$ to Caratheodory-measurable sets.

The dilation and translation operators on $L^2(\mathcal R,\mathfrak H^s)$ defined by 
$$Uf(x)=\frac{1}{\sqrt2}f\left(\frac x3\right),\quad Tf(x)=f(x-1),$$
are unitary and satisfy the covariance relation $UTU^{-1}=T^3$. Moreover $\varphi=\chi_{\mathbf C}$ is an orthogonal scaling function: it satisfies the scaling equation 
$$U\varphi=\frac{1}{\sqrt2}\left(\varphi+T^2\varphi\right),$$
its integer translates are orthogonal, and it generates a multiresolution, in the same sense as the one described above for $L^2(\br)$.

At the FL-IA-CO-OK Workshop in February 2009 in Iowa City, after discussions with Judy Packer and Palle Jorgensen, the following question arose: is this representation irreducible, i.e., is the commutant of $\{U,T\}$ trivial in $\mathcal B(L^2(\mathcal R,\mathfrak H^s))$?

This is one of the questions that motivated the investigation in the present paper. Even though we do not give a definite answer to this question, we will present some positive evidence that the respresentation is {\it not} irreducible. 

\subsection{Wavelet representations}
Although specific examples of wavelet representations have been studied for some time by many authors, a useful generalization
of this concept which can be used in a variety of situations was first introduced in \cite{DuJo07} to extend the multiresolution techniques to other discrete dynamical systems, and to construct orthonormal wavelet bases on other spaces beside $L^2(\br)$. The idea was to keep some of the essential properties of the multiresolutions mentioned above, but now as axioms in some abstract Hilbert space.

For more connections between wavelet representations, generalized multiresolutions and direct limits we refer to \cite{BCM02,BMM99,bfmp2,bfmp1,aijln,ijkln}.

Let $X$ be a compact metric space. Let $r:X\rightarrow X$ be a Borel measurable function and assume that $0<\#r^{-1}(x)<\infty$ for all $x\in X$. 
Assume that $\mu$ is a Borel probability measure on $X$ which is {\it strongly invariant}, i.e.,
\begin{equation}
\int f\,d\mu=\int\frac{1}{\#r^{-1}(x)}\sum_{r(y)=x}f(y)\,d\mu(x),\quad (f\in C(X)).
\label{eqid1}
\end{equation}
\end{definition}

\begin{theorem}\label{thi2}\cite[Corollary 3.6]{DuJo07}
Let $m_0$ be a function in $L^\infty(X,\mu)$ such that 
\begin{equation}
\frac{1}{\#r^{-1}(x)}\sum_{r(y)=x}|m_0(y)|^2=1,\quad(x\in X)
\label{eqit1}
\end{equation}
Then there exists a Hilbert space $\H$, a unitary operator $U$ on $\H$, a representation $\pi$ of $L^\infty(X)$ on $\H$ and an element $\varphi$ of $\H$ such that
\begin{enumerate}
\item 
(Covariance) $U\pi(f)U^{-1}=\pi(f\circ r)$ for all $f\in L^\infty(X)$. 
\item (Scaling equation) $U\varphi=\pi(m_0)\varphi$
\item (Orthogonality) $\ip{\pi(f)\varphi}{\varphi}=\int f\,d\mu$ for all $f\in L^\infty(X)$. 
\item (Density) $\{U^{-n}\pi(f)\varphi\,|\,n\in\bn, f\in L^\infty(X)\}$ is dense in $\H$.

\end{enumerate}
Moreover they are unique up to isomorphism.

\end{theorem}

\begin{definition}\label{defi3}
We say that $(\H,U,\pi,\varphi)$ in Theorem \ref{thi2} is the {\it wavelet representation} associated to $m_0$.
\end{definition}

The paper is structured as follows: in Section 2 we describe a concrete realization of the wavelet representation on the solenoid. This was mainly done in \cite{DuJo07}, but we present here a slightly different form. We show how the irreducibility of the wavelet representation is related to the ergodic properties of the shift on the solenoid, and to the fixed points of a transfer operator.

In Theorem \ref{thi6} we describe the multiresolution structure that comes with a wavelet representation.

In Section 3 we investigate two examples. The first one is the wavelet representation associated to an arbitrary map $r$, and the constant function $m_0=1$. Using the multiresolution structure we show in Theorem \ref{th7} that the shift on the solenoid is ergodic iff $r$ is ergodic.

The second example is the wavelet representation associated to the Cantor set, introduced in \cite{DuJo06}. That is $r(z)=z^3$ on the unit circle and $m_0(z)=\frac{1}{\sqrt{2}}(1+z^2)$. We show in Proposition \ref{pr2.1} that there is an $L^2(\bt,\mu)$ function which is a fixed point for the transfer operator $R_{m_0}$. However, this function is not bounded, and it does not satisfy the conditions of Theorem \ref{thi5.2}, so we cannot conclude that the representation is irreducible. In any case, this does provide some evidence that the representation might not be irreducible.

\section{Representations on the solenoid}

When the function $m_0$ is non-singular, i.e., $\mu(\{x\in X\,|\, m_0(x)=0\})=0$, the wavelet representation can be realized more concretely on the solenoid. We describe this realization. The basic idea is to regard the multiresolution as a martingale; the idea appeared initially in \cite{CoRa90} and \cite{Gun00}. It was then developed in \cite{DuJo07} for a larger class of maps $r$ and low-pass filters $m_0$ (see also \cite{Gun07}). Since we will need this representation in a slightly different form we include some of the details, and we refer to \cite{DuJo07} for a more rigurous account.

\begin{definition}\label{defi4}
Let 
\begin{equation}
X_\infty:=\left\{(x_0,x_1,\dots)\in X^{\bn}\,|\, r(x_{n+1})=x_n\mbox{ for all }n\geq 0\right\}
\label{eqi4_1}
\end{equation}
We call $X_\infty$ the {\it solenoid} associated to the map $r$.

On $X_\infty$ consider the $\sigma$-algebra generated by cylinder sets. 
Let $r_\infty:X_\infty\rightarrow X_\infty$
\begin{equation}
r_\infty(x_0,x_1,\dots)=(r(x_0),x_0,x_1,\dots)\mbox{ for all }(x_0,x_1,\dots)\in X_\infty
\label{eqi4_2}
\end{equation}
Then $r_\infty$ is a measurable automorphism on $X_\infty$. 

Define $\theta_0:X_\infty\rightarrow X$, 
\begin{equation}
\theta_0(x_0,x_1,\dots)=x_0.
\label{eqi4_3}
\end{equation}
The measure $\mu_\infty$ on $X_\infty$ will be defined by constructing some path measures $P_x$ on the fibers $\Omega_x:=\{(x_0,x_1,\dots)\in X_\infty\,|\, x_0=x\}.$

Let 
$$c(x):=\#r^{-1}(r(x)),\quad W(x)=|m_0(x)|^2/c(x),\quad(x\in X).$$
Then 
\begin{equation}
\sum_{r(y)=x}W(y)=1,\quad(x\in X)
\label{eqi4_4}
\end{equation}
$W(y)$ can be thought of as the transition probability from $x=r(y)$ to one of its pre-images $y$ under the map $r$.

For $x\in X$, the path measure $P_x$ on $\Omega_x$ is defined on cylinder sets by
\begin{equation}
P_x(\{(x_n)_{n\geq0}\in\Omega_x\,|\, x_1=z_1,\dots,x_n=z_n\})=W(z_1)\dots W(z_n)
\label{eqi4_5}
\end{equation}
for any $z_1,\dots, z_n\in X$.

This value can be interpreted as the probability of the random walk to go from $x$ to $z_n$ through the points $x_1,\dots,x_n$.

Next, define the measure $\mu_\infty$ on $X_\infty$ by 
\begin{equation}
\int f\,d\mu_\infty=\int_{X}\int_{\Omega_x}f(x,x_1,\dots)\,dP_{x}(x,x_1,\dots)\,d\mu(x)
\label{eqi4_7}
\end{equation}
for bounded measurable functions on $X_\infty$.

Consider now the Hilbert space $\H:=L^2(X_\infty,\mu_\infty)$. Define the operator 
\begin{equation}
Uf=m_0\circ\theta_0\,f\circ r_\infty,\quad(f\in L^2(X_\infty,\mu_\infty))
\label{eqi4_8}
\end{equation}

Define the representation of $L^\infty(X)$ on $\H$ 
\begin{equation}
\pi(f)g=f\circ\theta_0\,g,\quad(f\in L^\infty(X),g\in \H)
\label{eqi4_9}
\end{equation}

Let $\varphi=1$ be the constant function $1$.
\end{definition}

\begin{theorem}\label{thi5}
Suppose $m_0$ is non-singular, i.e., $\mu(\{x\in X\,|\, m_0(x)=0\})=0$. Then the data $(\H,U,\pi,\varphi)$ from Definition \ref{defi4} form the wavelet representation associated to $m_0$.
\end{theorem}

\begin{proof}
We check that $U$ is unitary, all the other relations follow from some easy computations. To check that $U$ is an isometry it is enough to apply it on functions $f$ on $X_\infty$ which depend only on the first $n+1$ coordinates $f=f(x_0,\dots,x_n)$. Then $f\circ r_\infty$ depends only on $x_0,\dots, x_{n-1}$. We have, using \eqref{eqi4_5} and the strong invariance of $\mu$:
$$\int |m_0\circ\theta_0|^2 |f\circ r_\infty|^2\,d\mu_\infty=$$$$
\int_X |m_0(x_0)|^2\sum_{r(x_1)=x_0,\dots, r(x_{n-1})=x_{n-2}}W(x_1)\dots W(x_{n-1})f(r(x_0),x_0,x_1,\dots,x_{n-1})\,d\mu(x_0)$$
$$=\int_X \frac{1}{\#r^{-1}(x)}\sum_{r(y)=x}|m_0(y)|^2\sum_{r(x_1)=y,r(x_2)=x_1,\dots, r(x_{n-1})=x_{n-2}}W(x_1)\dots W(x_{n-1})\cdot$$
$$\cdot f(r(y),y,x_1,\dots,x_{n-1})\,d\mu(x)=$$
$$\int_X\sum_{y_1,\dots y_n} W(y_1)\dots W(y_n)f(x,y_1,\dots,y_n)\,d\mu(x)=\int f\,d\mu_\infty.$$
This shows that $U$ is an isometry. 

The fact that $m_0$ is non-singular insures that $U$ is onto and has inverse 
$$Uf=\frac{1}{m_0\circ\theta_0\circ r_\infty^{-1}}f\circ r_\infty^{-1}$$
\end{proof}

 The commutant of the wavelet representations, i.e., the set of operators that commute with both the ``dilation'' operator $U$ and the ``translation'' operators $\pi(f)$,has a simple description that we will present below. Also the operators in the commutant are in one-to-one correspondence with bounded fixed points of the transfer operator. The commutant of the classical wavelet representation on $L^2(\br)$ was computed in \cite{DaLa98}. We will be interested in computing this commutant for other choices of filters, such as $m_0=1$ or for the wavelet representation associated to the Cantor set.

\begin{theorem}\cite[Theorem 7.2]{DuJo07}\label{thi5.1}
Suppose $m_0$ is non-singular and let $(\H,U,\pi,\varphi)$ be the wavelet representation as in Theorem \ref{thi5}. 
\begin{enumerate}
\item
The commutant $\{U,\pi\}'$ in $\mathcal B(\H)$ consists of operators of multiplication by functions $f\in L^\infty(X_\infty,\mu_\infty)$ which are invariant under $r_\infty$, i.e., $f\circ r_\infty=f$. We call these functions cocycles.
\item There is a one-to-one correspondence between cocycles and bounded fixed points for the transfer operator $R_{m_0}$ defined for functions on $X$:
\begin{equation}
R_{m_0}f(x)=\frac{1}{\#r^{-1}(x)}\sum_{r(y)=x}|m_0(y)|^2f(y),\quad(x\in X)
\label{eq5.1_1}
\end{equation}
The correspondence is defined as follows:

For a bounded cocycle $f$ on $X_\infty$ the function 
\begin{equation}
h(x)=\int_{\Omega_x}f(x,x_1,\dots)\,dP_x(x,x_1,x_2,\dots)
\label{eqi5.1_2}
\end{equation}
is a bounded fixed point for $R_{m_0}$, i.e., $R_{m_0}h=h$.

For a bounded measurable fixed point $h$ for the transfer operator $R_{m_0}$, the limit exists $\mu_\infty$-a.e.
\begin{equation}
f(x_0,x_1,\dots):=\lim_{n\rightarrow\infty}h(x_n),\quad((x_0,x_1,\dots)\in X_\infty)
\label{eq5.1_3}
\end{equation}
and defines a bounded cocyle.
\end{enumerate}
\end{theorem}

Next, we describe the multiresolution structure associated to a wavelet representation. The proof is standard in wavelet theory, but we include the main ideas for the benefit of the reader.
\begin{theorem}\label{thi6}
Let
$$V_0:=\Span\left\{\pi(f)\varphi\,|\,f\in L^\infty(X)\right\},$$
$$V_n:=U^{-n}V_0,\quad(n\in\bz).$$
Then
\begin{enumerate}
\item
$UV_0\subset V_0$.
\item 
$\cj{\bigcup_{n\in\bz}V_n}=\H.$
\item
$V_0$ is an invariant subspace for the representation $\pi$. The spectral measure of the representation $\pi$ restricted to $V_0$ is $\mu$ and the multiplicity function is constant 1. 
\item 
$V_1$ is an invariant subspace for the representation $\pi$. The spectral measure of the representation $\pi$ restricted to $V_1$ is $\mu$ and the multiplicity function is $\mathfrak m_{V_1}(x)=\#r^{-1}(x)$, $x\in X$.
\item 
Let $W_0:=V_1\ominus V_0$. Then $W_0$ is invariant for $\pi$. The multiplicity function of $\pi$ on $W_0$ is $\mathfrak m_{W_0}(x)=\#r^{-1}(x)-1$.
\item 
$$\left(\bigoplus_{n\in\bz}U^n W_0\right)\oplus \bigcap_{n\in\bz}V_n=\H.$$
\item 
Let $N:=\sup_{x\in X}\#r^{-1}(x)\in\bn\cup\{\infty\}$. There exists functions $\psi_1,\dots,\psi_N$ (if $N$ is $\infty$ then the functions $\psi$ are just indexed by natural numbers, we don't have a $\psi_\infty$) in $W_0$ with the following properties:
\begin{equation}
\ip{U^n\pi(f)\psi_i}{U^m\pi(g)\psi_j}=\delta_{mn}\delta_{ij}\int f\cj g \chi_{\{\#r^{-1}(x)\geq i+1\}}\,d\mu,\quad(f,g\in L^\infty(X),m,n\in\bz, i,j\in\{1,\dots, N\})
\label{eqi6_1}
\end{equation}
\begin{equation}
\Span\left\{U^n\pi(f)\psi_i\,|\,f\in L^\infty(X),n\in\bz, i\in\{1,\dots, N\}\right\}=\H\ominus\bigcap_{n\in\bz}V_n
\label{eqi6_2}
\end{equation}

\end{enumerate}
\end{theorem}

\begin{proof}

(i) follows from the scaling equation, (ii) follows from the desity property of the wavelet representation, (iii) follows from the orthogonality. 
The fact that $V_1$ is invariant for $\pi$ follows from the covariance relation. The multiplicity function for $V_1$ was computed in \cite[Theorem 4.1]{DuJo07}. (v) follows from (iv). (vi) follows from the fact that $U$ is unitary so $U^{-n}W_0=V_{n+1}\ominus V_n$ for all $n\in\bz$.

For (vii) consider the space
$$L^2(X,\mu,\mathfrak m_{W_0}):=\left\{f:X\rightarrow\cup_{x\in X}\bc^{\mathfrak m_{W_0}(x)}\,|\, f(x)\in \bc^{\mathfrak m_{W_0}(x)}\mbox{ for all }x\in X, \int_X \|f(x)\|^2\,d\mu(x)<\infty\right\}.$$
On this space we have the representation of $L^\infty(X)$ by multiplication $M_f$. By (v) there is an isomorphism $J:W_0\rightarrow L^2(X,\mu,\mathfrak m_{W_0})$ such that 
$J\pi(f)=M_fJ$ for all $f\in L^\infty(X)$.

Let $e_i$ be the canonical vectors in $\bc^n$. Define the functions $\eta_i\in L^2(X,\mu,\mathfrak m_{W_0})$:
$$\eta_i(x)=\left\{\begin{array}{cc}
e_i,&\mbox{ if } \mathfrak m_{W_0}(x)=\#r^{-1}(x)-1\geq i\\
0,&\mbox{ otherwise.}
\end{array}
\right.$$
Let $\psi_i:=J^{-1}\eta_i$. 

It is then easy to see that if $i\neq j$ then $\ip{\eta_i(x)}{\eta_j(x)}=0$ for all $x$, so $\ip{\pi(f)\psi_i}{\pi(g)\psi_j}=0$ for all $f,g\in L^\infty(X)$, $i\neq j$. Also
$$\ip{f\eta_i}{g\eta_i}=\int_{\{\#r^{-1}(x)-1\geq i\}}f\cj g\,d\mu.$$
This, together with (vi) implies \eqref{eqi6_1}.

Equation \eqref{eqi6_2} is also a consequence of (vi) if we show that $\pi(f)\psi_i$ span $W_0$. But it is clear that $M_f\eta_i$ span $L^2(X,\mu,\mathfrak m_{W_0})$ so, applying $J^{-1}$ we get the result.
\end{proof}

Finally, we present several equivalent formulations of the problem of the irreducibility of a wavelet representation.
\begin{theorem}\label{thi5.2}
Suppose $m_0$ is non-singular. The following affirmations are equivalent:
\begin{enumerate}
\item
The wavelet representation is irreducible, i.e., the commutant $\{U,\pi\}'$ is trivial.
\item
The automorphism $r_\infty$ on $(X_\infty,\mu_\infty)$ is ergodic.
\item
The only bounded measurable fixed points for the transfer operator $R_{m_0}$ are the constants.
\item
There does not exist a non-constant fixed point $h\in L^p(X,\mu)$ with $p>1$ of the transfer operator $R_{m_0}$ with the property that
\begin{equation}
\sup_{n\in\bn}\int_X|m_0^{(n)}(x)|^2|h(x)|^p\,d\mu(x)<\infty
\label{eqi5.2_1}
\end{equation}
where 
\begin{equation}
m_0^{(n)}(x)=m_0(x)m_0(r(x))\dots m_0(r^{n-1}(x)),\quad (x\in X).
\label{eq5.2_2}
\end{equation}
\item
If $\varphi'\in\H$, satisfies the same scaling equation as $\varphi$, i.e., $U\varphi'=\pi(m_0)\varphi'$, then $\varphi'$ is a constant multiple of $\varphi$.
\end{enumerate}

\end{theorem}

\begin{proof}
The equivalences of (i)--(iii) follow immediately from Theorem  \ref{thi5.1}. It is also clear that (iv) implies (iii), because bounded functions satisfy \eqref{eqi5.2_1} with any $p>1$. Indeed, using the strong invariance of $\mu$:
$$\int |m_0^{(n)}|^2|h|^p\,d\mu\leq \|h\|_\infty\int |m_0^{(n)}|^2\,d\mu=\|h\|_\infty\int_X R_{m_0}^n1\,d\mu=\|h\|_\infty.$$  

We prove that (ii) implies (iv) by contradiction. Suppose there is a non-constant $h$ with the given properties. Define the functions on $X_\infty$
$$h_n(x_0,x_1,\dots)=h(x_n),\quad(x_0,x_1,\dots)\in X_\infty.$$
Then $(h_n)_n$ is a martingale with respect to the filtration $\theta_n^{-1}(\mathcal B)$, where $\mathcal B$ is the Borel $\sigma$-algebra in $X$ and $\theta_n:X_\infty\rightarrow X$, $\theta_n(x_0,x_1,\dots)=x_n$. We denote by $\be_n$ the conditional expectation onto $\theta_n^{-1}(\mathcal B)$. We have, since $h_{n+1}$ depends only on $x_0,\dots,x_{n+1}$:
$$\be_{n}(h_{n+1})(x_0,\dots,x_n,\dots)=\frac{1}{\#r^{-1}(x_n)}\sum_{r(x_{n+1})=x_n}|m_0(x_{n+1})|^2h_{n+1}(x_0,\dots,x_{n+1},\dots)=$$$$ 
\frac{1}{\#r^{-1}(x_n)}\sum_{r(x_{n+1})=x_n}|m_0(x_{n+1})|^2h(x_{n+1})=h(x_n)=h_n(x_0,x_1,\dots).$$

We want to apply Doob's discrete martingale convergence theorem. We have to check that 
\begin{equation}
\sup_n\int_{X_\infty} |h_n|^p\,d\mu_\infty<\infty.
\label{eqi5.2_2}
\end{equation}
But, using the strong invariance of $\mu$ applied $n$ times:
$$\int_{X_\infty} |h_n|^p\,d\mu_\infty=\int_X\sum_{r(x_1)=x_0,\dots r(x_n)=x_{n-1}}W(x_1)\dots W(x_n)|h(x_n)|^p\,d\mu(x_0)=$$
$$=\int_X R_{m_0}^n|h|^p\,d\mu=\int_X |m_0^{(n)}|^2|h|^p\,d\mu$$
Doob's theorem implies then that $$f(x_0,x_1,\dots)=\lim_n h_n(x_0,x_1,\dots)$$ exists $\mu_\infty$-a.e., and in $L^1(X_\infty,\mu_\infty)$. Then 
$$\be_0(f)=\lim_n \be_0(h_n)=h$$
so $f$ is not a constant. But we also have 
$$f\circ r_\infty(x_0,x_1,\dots)=f(r(x_0),x_0,x_1,\dots)=\lim_n h(x_{n-1})=f(x_0,x_1,\dots)$$
$\mu_\infty$-a.e. This contradicts the fact that $r_\infty$ is ergodic.

$(ii)\Rightarrow(v)$. Take a $\varphi'$ as in (v). Then, the scaling equation implies $$m_0\circ\theta_0\,\varphi'\circ r_\infty=U\varphi'=\pi(m_0)\varphi'=m_0\circ\varphi'.$$ Since $m_0$ is non-singular, this implies that $\varphi'\circ r_\infty=\varphi'$. But since $r_\infty$ is ergodic it follows that $\varphi'$ is a constant, i.e., $\varphi'$ is a constant multiple of $\varphi$.

$(v)\Rightarrow(ii)$. If $r_\infty$ is not ergodic, then one can take $\varphi'$ to be the characteristic function of a proper $r_\infty$-invariant set. It follows immediately that $\varphi'$ satisfies the scaling equation, and thus its existence contradicts (v).
\end{proof}

\section{Examples}

In this section we will consider two examples. The first example is the wavelet representation associated to $m_0=1$. The map $r$ can be any map satisfying the conditions above. We show that the wavelet representation associated to $m_0=1$ is irreducible if and only if $r$ is ergodic. 

The second example is the wavelet representation associated to the Cantor set, representation that was defined in \cite{DuJo06b}. The representation is associated to the map $r(z)=z^3$, for $z\in\bc$, $|z|=1$, and the QMF filter $m_0(z):=(1+z^2)/\sqrt2$. While we were not able to determine if this representation is irreducible or not, we present several equivalent formulations of the problem, in terms of the existence of solutions for refinement equations or the existence of fixed points for transfer operators. We find a non-trivial fixed point for the associated tranfer operator which is in $L^2(\bt)$, but it is not bounded (so it does not settle the problem, but gives some positive evidence that the representation might be reducible). At the same time we show that it is hard to give a constructive solution for the irreducibility problem: in Proposition \ref{prsca} we prove that the refinement equation has no non-trivial compactly supported solutions. In Corollary \ref{conol1} we show that the transfer operator has no non-trivial solutions with Fourier transform in $l^1(\bz)$. In Proposition \ref{prcasc} we show that the method of successive approximations will not produce a new solution to the refinement equation, if the seed is compactly supported.

\subsection{The wavelet representation associated to $m_0=1$}
\begin{theorem}\label{th7}
Let $m_0=1$ and let $(\H,U,\pi,\varphi)$ be the associated wavelet representation. The following affirmations are equivalent:
\begin{enumerate}
\item The automorphism $r_\infty$ on $(X_\infty,\mu_\infty)$ is ergodic.
\item The wavelet representation is irreducible.
\item The only bounded functions which are fixed points for the transfer operator $R_1$, i.e.,
$$R_1h(x):=\frac{1}{\#r^{-1}(x)}\sum_{r(y)=x}h(y)=h(x)$$
are the constant functions.
\item The only $L^2(X,\mu)$-functions which are fixed points for the transfer operator $R_1$, are the constants. 
\item The endomorphism $r$ on $(X,\mu)$ is ergodic. 
\end{enumerate}

\end{theorem}

\begin{proof}
The equivalence of (i)--(iv) is given in Theorem \ref{thi5.2}. We will prove that (i) and (iv) are equivalent. 

$(i)\Rightarrow(v)$. Suppose $r$ is not ergodic. Let $f$ be a bounded, non-constant $\mu$-a.e., function on $X$ such that $f=f\circ r$. Define $\tilde f:=f\circ\theta_0$. 
Then it is easy to see that $\tilde f=\tilde f\circ r_\infty$. But since $r_\infty$ is ergodic this implies that $\tilde f$ is constant $\mu_\infty$-a.e.  But since $\tilde f=f\circ\theta_0$ depends only on the first coordinate, this implies that $f$ is constant $\mu$-a.e.

$(v)\Rightarrow(i)$. Let $f$ be a bounded function on $X_\infty$ such that $f=f\circ r_\infty$. We use Theorem \ref{thi6}. Pick $g\in L^\infty(X)$, and $i\in\{1,\dots,N\}$ arbitrary. Assuming that $\pi(g)\psi_i\neq 0$, let $A:=\|\pi(g)\psi_i\|$. (The case $A=0$ can be treated easily) Then we see that for all $n\in\bz$ we have 
$$\ip{f}{U^n\frac{1}{A}\pi(g)\psi_i}=\ip{U^{-n}f}{\frac{1}{A}\pi(g)\psi_i}=\ip{f\circ r_\infty^{-n}}{\frac{1}{A}\pi(g)\psi_i}
=\ip{f}{\frac{1}{A}\pi(g)\psi_i}$$
Thus these numbers do not depend on $n$. Moreover, we know that as $n$ varies, the vectors $U^n\frac{1}{A}\pi(g)\psi_i$ are orthogonal. Using Bessel's inequality, we have 
$$\infty\cdot\left|\ip{f}{\frac{1}{A}\pi(g)\psi_i}\right|^2=\sum_{n\in\bz}\left|\ip{f}{U^n\frac{1}{A}\pi(g)\psi_i}\right|^2\leq\|f\|^2<\infty.$$
This implies that all these numbers $\ip{f}{U^n\frac{1}{A}\pi(g)\psi_i}$ have to be 0.

Thus $f$ is orthogonal to all $U^n\pi(g)\psi_i$, and, by Theorem \ref{thi6}(vii), this shows that $f\in \cap_n V_n$. In particular $f\in V_0$ so there exists a function $\tilde f\in L^2(X,\mu)$ such that $f=\tilde f\circ\theta_0$. But since $f$ is invariant under $r_\infty$, $\tilde f$ is invariant under $r$ so it has to be constant $\mu$-a.e., so $f$ is constant $\mu_\infty$-a.e. Therefore $\mu_\infty$ is ergodic.
\end{proof}

\subsection{The wavelet representation associated to the Cantor set}
Recall (\cite{DuJo06b}) that the wavelet representation associated to the Cantor set is associated to $r(z)=z^3$ on the unit circle $\bt$, and the function 
\begin{equation}
m_0(z)=\frac{1}{\sqrt{2}}(1+z^2),\quad(z\in\bt)
\label{eqc1}
\end{equation}
As we mentioned in the introduction, in section \ref{secwc}, it can be realized on the Hilbert space $L^2(\mathcal R,\mathfrak H^s)$ associated to the Hausdorff measure $\mathfrak H^s$ on the subset $\mathcal R$.

\begin{theorem}\label{thirr}
The following assertions are equivalent:
\begin{enumerate}
\item
The wavelet representation associated to $m_0$ is irreducible. 
\item 
If a sequence $(a_k)_{k\in\bz}\in l^2(\bz)$ satisfies the properties that $\sum_{k\in\bz}a_kz^k\in L^\infty(\bt,\mu)$ and 
\begin{equation}
a_k=\frac12a_{3k-2}+a_{3k}+\frac12 a_{3k+2},\quad(k\in\bz)
\label{eqirr2}
\end{equation}
then $a_k=0$ for all $k\neq 0$.
\item 
If a function $\xi\in L^2(\mathcal R,\mathfrak H^s)$ satisfies the refinement equation
$$\xi(x)=\xi(3x)+\xi(3x-2),\mbox{ for $\mathfrak H^s$-a.e. $x\in\mathcal R$},$$
then $\xi$ is a constant multiple of the characteristic function of the Cantor set $\mathbf C$.

\end{enumerate}

\end{theorem}

\begin{proof}
To prove (i)$\Leftrightarrow$(ii) we use the equivalence of  (i) and (ii) in Theorem \ref{thi5.2} and the following Lemma.
\begin{lemma}\label{lemfpc}
Let $f\in L^2(\bt,\mu)$, $f=\sum_{k\in\bz}f_kz^k$. Then $f$ is a fixed point for the transfer operator $R_{m_0}$ iff
\begin{equation}
f_n=\frac12 f_{3n-2}+ f_{3n}+\frac12 f_{3n+2},\quad(n\in\bz)
\label{eqc4}
\end{equation}
\end{lemma}

\begin{proof}
We have
\begin{equation}
|m_0(z)|^2=1+\frac12z^2+\frac12 z^{-2}
\label{eqc2}
\end{equation}

Using the strong invariance of $\mu$, we compute the Fourier coefficients of $R_{m_0}f$ for a function $f\in L^2(\bt,\mu)$:
$$(R_{m_0}f)_k=\ip{R_{m_0}f}{z^k}=\int_{\bt} R_{m_0}f\cdot z^{-k}\,d\mu=\int_{\bt} \frac{1}{3}\sum_{w^3=z}|m_0(w)|^2f(w)\cdot w^{-3k}\,d\mu(z)$$
$$=\int_{\bt}|m_0(z)|^2f(z)z^{-3k}\,d\mu(z)=\int_{\bt}\left(z^{-3k}+\frac12 z^{-(3k-2)}+\frac12 z^{-(3k+2)}\right) f(z)\,d\mu(z)=\frac12 f_{3k-2}+ f_{3k}+\frac12 f_{3k+2}$$
Thus

\begin{equation}
(R_{m_0}f)_k=\frac12 f_{3k-2}+ f_{3k}+\frac12 f_{3k+2},\quad(k\in\bz)
\label{eqc3}
\end{equation}
This implies \eqref{eqc4}
\end{proof}

To see that (i) and (iii) are equivalent, use (v) in Theorem \ref{thi5.2}. 

\end{proof}

Next, we will analyze conditions (ii) and (iii) in Theorem \ref{thirr} and rule out some solutions. More precisely, in Propositon \ref{prsca} we prove that there are no compactly supported solutions for the refinement equation in (iii); in Corollary \ref{conol1} we show that there are no $l^1$-solutions for the fixed point problem in (ii). However, in Proposition \ref{pr2.1} we do find an $l^2$-solution. In Proposition \ref{prcasc} we show that the method of successive approximations produces highly divergent sequences for the refinement equation in (iii).

\begin{proposition}\label{prsca}
The only  Borel measurable solutions for the refinement equation
$$\varphi(x)=\varphi(3x)+\varphi(3x-2),\quad(x\in\mathcal R)$$
with bounded support, are constant multiples of the characteristic function of the Cantor set $\mathbf C$, up to $\mathfrak H^s$-measure zero.
\end{proposition}

\begin{proof}
Let $a:=\sup\{x\in\mathcal R\,|\,\varphi(x)\neq0\}$. We cannot have $a>1$, because then there exists a sequence $x_n\leq a$ that converges to $a$ and such that $\varphi(x_n)\neq 0$. But then either $\varphi(3x_n)$ or $\varphi(3x_n-2)$ is non-zero, and both $3x_n$ and $3x_n-2$ are bigger than $a$ for $n$ large. Thus $a\leq 1$. 
A similar argument shows that 0 is a lower bound for the support of $\varphi$. Thus $\varphi$ has to be supported on $[0,1]$. Let $K$ be its support, i.e.,
$K$ is the closure in $\br$ of $\{x\in\mathcal R\,|\,\varphi(x)\neq 0\}$. We claim that 
\begin{equation}
K=\frac{K}{3}\cup\frac{K+2}3
\label{eqsca1}
\end{equation}

If $x\in[0,1]$ and $\varphi(x)\neq0$ then either $\varphi(3x)$ or $\varphi(3x-2)$ is non-zero, therefore either $x\in K/3$ or $x\in (K+2)/3$. This proves one inclusion.

From the scaling equation, we have that 
$$\varphi(x/3)=\varphi(x)+\varphi(x-2)$$
But if $x\in[0,1]$, then $x-2$ is not, so 
$\varphi(x/3)=\varphi(x)$ for $x\in[0,1]$. 
Similarly $\varphi((x+2)/3)=\varphi(x)$ for $x\in[0,1]$. 

If $x\in K/3$ then $\varphi(3x)\neq 0$ and $3x\in[0,1]$, so $\varphi(x)=\varphi(3x)\neq0$, so $x\in K$. Hence $K/3\subset K$. Similarly $(K+2)/3\subset K$. This proves \eqref{eqsca1}. Since the Cantor set $\mathbf C$ is the only compact solution for \eqref{eqsca1} (see e.g. \cite{Hut81}), it follows that $\varphi$ is supported on the Cantor set.

The map $r(x)= 3x\mod 1$ on the Cantor set with the Hausdorff measure $\mathfrak H^s$, is ergodic, since it is conjugate to the shift on the symbolic space $\{0,1\}^\bn$, 
$\sigma(d_1,d_2,\dots)=(d_2,d_3,\dots)$ with the product measure, where $0$ and $1$ get equal probabilities $1/2$. The conjugating map is $\Psi(d_1,d_2,\dots)=\sum_{n\geq1}2d_n/3^n$.

Moreover $\varphi$ is invariant under the shift since $\varphi(x/3)=\varphi((x+2)/3)=\varphi(x)$ for $x\in \mathbf C$. Then, $\varphi$ must be constant on $\mathbf C$, and the proposition is proved.

\end{proof}

To study solutions for the fixed-point problem in Theorem \ref{thi5.2}(iii) or its particular form in Theorem \ref{thirr} (ii), we need some background on the transfer operator. The next theorem is contained in \cite{DuJo06b}, Theorem 5.1, Proposition 7.1, and Theorem 7.4.

\begin{theorem}\label{thruelle}\cite{DuJo06b}
Let $m_0(z)=\frac{1+z^2}{\sqrt{2}}$ and let $R_{m_0}$ be the corresponding transfer operator. 
\begin{enumerate}
\item
If $h\in C(\bt)$ and $R_{m_0}h=h$ then $h$ is constant.
\item
There are no functions $f\in C(\bt)$ and $\lambda\in\bc$ with $|\lambda|=1$, $\lambda\neq 1$ and $R_{m_0}f=\lambda f$. 
\item There is a unique Borel probability measure on $\bt$ such that 
$$\int_{\bt}R_{m_0}f\,d\nu=\int_{\bt}f\,d\nu,\quad(f\in C(\bt)).$$
Moreover $\nu$ has full support, in other words, every non-empty open subset of $\bt$ has positive measure. 
\item 
For all $f\in C(\bt)$, $\lim_{n\rightarrow\infty}R_{m_0}^nf=\nu(f)$, uniformly on $\bt$.
\end{enumerate}

\end{theorem}

\begin{corollary}\label{conol1}
There is no non-trivial solution for equation \eqref{eqirr2} in $l^1(\bz)$. By trivial, we mean a sequence $(a_k)_{k\in\bz}$ with $a_k=0$ for all $k\neq 0$.
\end{corollary}

\begin{proof}
Suppose $(a_k)_{k\in\bz}$ is a solution for \eqref{eqirr2} in $l^1$. Then $\sum_{k\in\bz}a_kz^k$ is uniformly convergent to a continuous function $h$, and $R_{m_0}h=h$. Then, by Theorem \ref{thruelle}, it follows that $h$ is a constant, so the sequence $(a_k)_{k\in\bz}$ is the trivial solution.

\end{proof}

In the next proposition we present a solution in $l^2(\bz)$ for equation \eqref{eqirr2}. However, its Fourier transform, while in $L^2(\bt,\mu)$, is not bounded, and therefore it does not offer a solution to our problem. It just gives some evidence that this wavelet representation might not be irreducible.

\begin{proposition}\label{pr2.1}
Define the sequence $(a_n)_{n\in\bz}$ as follows:

\begin{equation}
a_n:=\left\{\begin{array}{cc}
\frac1{2^k},&\mbox{ if $n$ is an even number between $3^k+1$ and $3^{k+1}-1$, $k\geq0$}\\
-\frac1{2^k},&\mbox{ if $n$ is an even number between $-(3^{k+1}-1)$ and $-(3^{k}+1)$, $k\geq0$}\\
0,&\mbox{ otherwise}.
\end{array}
\right.
\label{eqc5}
\end{equation}
Then the function 
\begin{equation}
h(z):=\sum_{k\in\bz}a_kz^k,\quad(z\in\bt)
\label{eq2.1_1}
\end{equation}
satisfies the following properties:
\begin{enumerate}
\item
$h\in L^2(\bt,\mu)$ but $h\not\in L^\infty(\bt,\mu)$.
\item $R_{m_0}h=h$.
\item
$$\sup_{n}\int_{\bt}|m_0^{(n)}|^2|h|^2\,d\mu=\infty.$$

\end{enumerate}
\end{proposition}
\begin{proof}
First we claim that $(a_n)_{n\in\bz}$ is in $l^2(\bz)$. Indeed, there are $3^k$ even numbers between $3^k+1$ and $3^{k+1}-1$. Then 
$$\sum_{n\in\bz}|a_n|^2=2\cdot\sum_{k\geq0}\left(\frac{1}{2^k}\right)^2\cdot 3^k=2\cdot\sum_{k=0}^\infty\left(\frac34\right)^k<\infty.$$
Thus $h\in L^2(\bt,\mu)$.

Next, we check that $R_{m_0}h=h$. Using Lemma \ref{lemfpc} we have to check that $(a_n)_{n\in\bz}$ satisfies equation \eqref{eqc4}. If $n$ is odd, then $3n,3n-2,3n+2$ are all odd, so the equation holds. 
If $n$ is even we have three cases. If $n=0$ then $a_{-2}=-1$, $a_2=1$, and the equation holds. Assume now $n$ is even and $n>0$. If $n$ is between $3^k+1$ and $3^{k+1}-1$. Then $3n-2$ is bigger than $3^{k+1}+1$ and $3n+2$ is less than $3^{k+2}-1$. And of course $3n,3n+2,3n-2$ are all even. Since we have 
$$a_n=\frac1{2^k}, a_{3n-2}=a_{3n}=a_{3n+2}=\frac{1}{2^{k+1}}$$
we see that the equation \eqref{eqc4} holds. 

The case $n<0$ can be treated similarly.

To prove (iii), we estimate the integral in \eqref{eq2.1_1}. This is the square of the $L^2$-norm of the function $f^{(n)}:=m_0^{(n)}h$, which can be computed as the sum of the squares of its Fourier coefficients, which we denote by $(a_k^{(n)})_{k\in\bz}$. 

We have $a_k^{(0)}=a_k$ for all $k$. Also,
$f^{(n+1)}=m_0(z^{3^{n}})f^{(n)}$ so 
\begin{equation}
a_k^{(n+1)}=\frac{a_k^{(n)}+a_{k-2\cdot 3^n}^{(n)}}{\sqrt{2}}.
\label{eq2.1_3}
\end{equation}

We prove by induction, that for all $n\geq 0$, and all $k\geq 3^n$, $k$ even, the sequence $(a^{(n)}_k)_k$ is decresing and non-negative.
For $n=0$, this is clear. Assume this holds for $n$ and prove it for $n+1$.  We have for $k$ even, and $k\geq 3^{n+1}$, $k-2\cdot 3^n\geq 3^n$ and is even. Then
$$a_{k+2}^{(n+1)}=\frac{a_{k+2}^{(n)}+a_{k+2-2\cdot 3^n}^{(n)}}{\sqrt{2}}\leq \frac{a_k^{(n)}+a_{k-2\cdot 3^n}^{(n)}}{\sqrt{2}}=a_k^{(n+1)}$$
and from the formula \eqref{eq2.1_3} it is clear that $a_k^{(n)}\geq 0$.

Next, we claim that for $k\geq 3^n$, even, 
\begin{equation}
a_k^{(n)}\geq \sqrt{2}^n a_k.
\label{eq2.1_4}
\end{equation}
Indeed, since $k-2\cdot 3^{n-1}\geq 3^{n-1}$, and $a^{(n-1)}$ is decreasing:
$$a_k^{(n)}=\frac{a_k^{(n-1)}+a_{k-2\cdot 3^{n-1}}^{(n-1)}}{\sqrt{2}}\geq \frac{2a_k^{(n-1)}}{\sqrt{2}}=\sqrt{2}a_k^{(n-1)}$$
Then, by induction $a_{k}^{(n)}\geq \sqrt{2}^na_k^{(0)}=\sqrt2^na_k$ for $k\geq 3^n$ even.

Now, using \eqref{eq2.1_4}, we have 
$$\|f^{(n)}\|^2=\sum_{k\in\bz}|a_k^{(n)}|^2\geq \sum_{k\geq 3^n}|a_k^{(n)}|^2\geq 2^n\sum_{k\geq 3^n}|a_k|^2=$$
$$2^n\sum_{m\geq n}\sum_{3^m\leq k< 3^{m+1}}|a_k|^2=2^n\sum_{m\geq n}3^m\left(\frac1{2^m}\right)^2=2^n\left(\frac34\right)^n\cdot\frac1{1-3/4}\rightarrow\infty$$
This proves (iii).

(iii) also implies that $h$ cannot be bounded, otherwise, using the strong invariance of $\mu$: 
$$\int_\bt |m_0^{(n)}|^2|h|^2\,d\mu\leq \|h\|_\infty\int_\bt R_{m_0}^n1\,d\mu=\|h\|_\infty.$$

\end{proof}

\begin{remark}
We know that the operators $U$ and $T$ satisfy the commutation relation $UTU^{-1}=T^3$. this implies that a formal series $\sum_{k\in\bz} T^{3^k}$ commutes with both $U$ and $T$. The problem with this series is that it is pointwise divergent at many points. For example, if $f$ has bounded support then the functions $T^{3^k}f$ will be disjointly supported for $k$ big enough, but will have the same $L^2(\mathcal R,\mathfrak H^s)$-norm, since $T$ is unitary. However, it is possible that the geometry of the space $L^2(\mathcal R,\mathfrak H^s)$ allows this formal series to be convergent on a large subspace, in which case an application of the spectral theorem for unbounded operators might prove that the representation is in fact not irreducible.

This remark and the existence of fixed points for the transfer operator in Proposition \ref{pr2.1} give us some positive evidence that the wavelet representation associated to the Cantor set is not irreducible. On the other hand Proposition \ref{prsca}, Corollary \ref{conol1} and the next Proposition \ref{prcasc} show that a constuctive solution will be hard to come by.

\end{remark}

 One way to try to obtain solutions for the refinement equation is to iterate the cascade operator.
\begin{definition}
The operator $M:=U^{-1}\pi(m_0)$ on $\H$ is called the cascade operator. 
\end{definition}
We prove that convergence of the iterates of the cascade cannot be obtained if one starts with a function with bounded support. 

\begin{proposition}\label{prcasc}
Let $\xi\in L^2(\mathcal R,\mathfrak H^s)$ with bounded support. Suppose $\xi$ is not a constant multiple of $\chi_{\mathbf C}$. Then there is a positive constant $c_\xi>0$ such that 
$$\lim_{n\rightarrow\infty}\|M^{n+1}\xi-M^n\xi\|^2=c_\xi.$$
In particular, the sequence $(M^n\xi)_{n\in\bn}$ is not convergent. 
\end{proposition}

\begin{proof}
First, we need to introduce the {\it correlation} function for $\xi_1,\xi_2\in \H$. This is defined by considering the representation on the solenoid. 
\begin{equation}\label{eqcor}
p(\xi_1,\xi_2)(x):=\int_{\Omega_x}\xi(x,x_1,\dots)\cj{\xi}_2(x,x_1,\dots)\,dP_x(x,x_1,\dots),\quad(x\in\bt).
\end{equation}
Note that the correlation function is in $L^1(\bt,\mu)$ and has the following property (and it is completely determined by it):
\begin{equation}\label{eqcor2}
\ip{\pi(f)\xi_2}{\xi_2}=\int_{\bt}fp(\xi_1,\xi_2)\,d\mu,\quad(f\in L^\infty(X,\mu))
\end{equation}

Moreover, we claim that 
\begin{equation}\label{eqcor3}
p(M\xi_1,M\xi_2)=R_{m_0}p(\xi_1,\xi_2)
\end{equation}
Indeed, we have
$$\int_{\bt} fp(M\xi_1,M\xi_2)\,d\mu=\ip{\pi(f)M\xi_1}{M\xi_2}=\ip{\pi(|m_0|^2f\circ r)\xi_1}{\xi_2}=$$$$\int_{\bt}|m_0|^2f\circ rp(\xi_1,\xi_2)\,d\mu=\int_{\bt}fR_{m_0}p(\xi_1,\xi_2)\,d\mu.$$

Now, take $\xi\in \H$ with bounded support, and not a constant multiple of $\chi_{\mathbf C}$. Then $M\xi-\xi$ is also of bounded support. 
We have
\begin{equation}
\|M^{n+1}\xi-M^n\xi\|^2=\int_{\bt}p(M^{n+1}\xi-M^n\xi,M^{n+1}\xi-M^n\xi)\,d\mu=\int_{\bt}R_{m_0}^n p(M\xi-\xi,M\xi-\xi)\,d\mu.
\label{eqcor4}
\end{equation}

If $\eta\in \H$ is a function of bounded support then, by \eqref{eqcor2}, we have that 
$$\int_{\bt}z^kp(\eta,\eta)\,d\mu=\ip{T^k\eta}{\eta},\quad(k\in\bz).$$
Therefore $p(\eta,\eta)\geq0$ is a trigonometric polynomial.

Thus $h_0:=p(M\xi-\xi,M\xi-\xi)\geq0$ is a trigonometric polynomial. We claim first that $h_0$ cannot be identically $0$. If that is the case then from \eqref{eqcor2} it follows that $\|M\xi-\xi\|^2=0$ so $M\xi=\xi$. But we saw in Proposition \ref{prsca} that the only solutions of the refinement equation that have bounded support are multiples of $\chi_{\mathbf C}$.

Since $h_0$ is not identically zero and $h_0\geq0$ and it is continuous it follows that $\nu(h_0)>0$, since $\nu$ has full support by Theorem \ref{thruelle}.
From \eqref{eqcor4}, using Theorem \ref{thruelle} and the fact that $h_0$ is continuous it follows that 
$\|M^{n+1}\xi-M^n\xi\|^2\rightarrow \int_{\bt}\nu(h_0)\,d\mu=\nu(h_0)>0$, and the result is obtained.
\end{proof}

\begin{remark}
In the interval of time between the submission and the acceptance of this paper, the first and third author have proved that the wavelet representation associated to the middle-third Cantor set is actually reducible \cite{DS10}. The proof is not constructive, so it is not clear how the operators in the commutant, or the $L^\infty$-fixed points of the transfer operator look like. The present paper shows that a constructive approach can be quite complicated. 
\end{remark}

\begin{acknowledgements}
We would like to thank professors Palle Jorgensen and Judith Packer for discussions and suggestions that motivated this paper. 
\end{acknowledgements}
\newcommand{\etalchar}[1]{$^{#1}$}
\def\cprime{$'$}

\end{document}